\documentclass[a4paper,12pt]{amsart}
\usepackage{amssymb,amsthm}

\input xy
\xyoption{all}

\begin{document}

\newtheoremstyle{Def}{8pt}{8pt}{}{}{\bfseries}{.}{.5em}{\thmname{#1}\thmnote{ #3}}

\newcommand{\D}{\ensuremath{\mathcal{D}}}
\newcommand{\C}{\mathbb{C}}
\renewcommand{\P}{\mathbb{P}}
\newcommand{\SL}{\operatorname{SL}}
\renewcommand{\H}{\operatorname{H}}
\newcommand{\SHom}{\mathcal{H}om}
\newcommand{\F}{\ensuremath{\mathcal{F}}}
\newcommand{\SH}{\mathcal{H}}
\newcommand{\LW}{\mathcal{L}(v)}
\newcommand{\g}{\mathfrak{g}}
\newcommand{\height}{\operatorname{ht}}
\newcommand{\kar}{\operatorname{char}}
\newcommand{\Spec}{\operatorname{Spec}}
\newcommand{\Supp}{\operatorname{Supp}}
\newcommand{\Sing}{\operatorname{Sing}}
\newcommand{\Gr}{\operatorname{Gr}}
\newcommand{\cd}{\operatorname{cd}}

\newcommand{\Fbb}{\mathbb{F}}
\newcommand{\Gbb}{\mathbb{G}}
\newcommand{\Qbb}{\mathbb{Q}}
\newcommand{\Rbb}{\mathbb{R}}
\newcommand{\Zbb}{\mathbb{Z}}
\newcommand{\Lcal}{\mathcal L}
\newcommand{\Mcal}{\mathcal M}
\newcommand{\Ocal}{\mathcal O}
\newcommand{\diag}{\operatorname{diag}}
\newcommand{\GL}{\operatorname{GL}}

\newtheorem{Corollary}{Corollary}[section]
\newtheorem{Lemma}{Lemma}[section]
\newtheorem{Theorem}{Theorem}[section]
\newtheorem{Proposition}[Theorem]{Proposition}
{\theoremstyle{Def} \newtheorem{Definition}{Definition}}
\newtheorem{Remark}{Remark}
\newenvironment{Proof}{{\it Proof.\/}}{\hfill $\square$\medskip}

\title[$F$-regularity of large Schubert varieties]
{$F$-regularity of large Schubert varieties} 
\author{Michel Brion \and Jesper Funch Thomsen}
\address{Institut Fourier, B.P.74\\
38402 Saint-Martin d'H\`eres Cedex, France}
\email{Michel.Brion@ujf-grenoble.fr}
\address{Institut for matematiske fag\\Aarhus Universitet\\ \AA rhus, Denmark}
\email{funch@imf.au.dk}

\begin{abstract}
Let $G$ denote a connected reductive algebraic group over an
algebraically closed field $k$ and let $X$ denote a projective 
$G\times G$-equivariant embedding of $G$. The large Schubert varieties
in $X$ are the closures of the double cosets $B g B$, where $B$
denotes a Borel subgroup of $G$, and $g \in G$. We prove that these
varieties are globally $F$-regular in positive characteristic,
resp.~of globally $F$-regular type in characteristic $0$. As a
consequence, the large Schubert varieties are normal and
Cohen-Macaulay.
\end{abstract}

\maketitle

\section{Introduction}

The class of globally $F$-regular varieties was introduced by Smith in
\cite{Smith} ; these are projective algebraic varieties in positive
characteristics such that all the ideals in their homogeneous
coordinate rings are tightly closed. The globally $F$-regular
varieties (and their analogues in characteristic $0$, the varieties of
globally $F$-regular type) have remarkable properties, e.g., they are
normal and Cohen-Macaulay, and the higher cohomology groups of all nef
invertible sheaves are trivial. 

Examples of globally $F$-regular varieties include the projective
toric varieties (Prop.6.4 in \cite{Smith}). In this note, we obtain
the global $F$-regularity of a wider class of varieties with
algebraic group action: the $G\times G$-equivariant projective
embeddings of any connected reductive group $G$, and the closures in
any such embedding of the double cosets $B g B$, where $B$ denotes a
Borel subgroup of $G$, and $g\in G$ is arbitrary. By the Bruhat
decomposition, these ``large Schubert varieties'' are parametrized by
the Weyl group of $G$; examples are the closures of parabolic
subgroups. We also show that the large Schubert varieties are of
globally $F$-regular type in characteristic $0$.

For this, we exploit the close relation between global $F$-regularity
and Frobenius splitting established in \cite{Smith}, and the  
Frobenius splitting properties of the large Schubert varieties,
proved in \cite{BriPol} and \cite{Rit2}. Another key ingredient is the 
global $F$-regularity of the flag varieties and their Schubert
varieties (\cite{LauPedTho}). Note that, unlike for Schubert
varieties, no desingularization of large Schubert varieties is known
in general; this makes our arguments somewhat indirect. 

As a consequence of our result, the large Schubert varieties in any 
equivariant embedding of $G$ are normal and Cohen-Macaulay. This
was first proved in the case of the canonical compactification of a 
semisimple adjoint group, by Frobenius splitting methods (\cite{BriPol}). 
Then Rittatore showed that all the equivariant embeddings of connected
reductive groups are Cohen-Macaulay, again by Frobenius splitting
methods (\cite{Rit2}). On the other hand, the Cohen-Macaulayness of
large Schubert varieties in the space of $n\times n$ matrices
(regarded as an equivariant embedding of the general linear group
$\GL_n$) was established by Knutson and Miller via a degeneration
argument, see \cite{KnuMil}. 

The group $G$, regarded as a homogeneous space under $G\times G$, is
an example of a spherical homogeneous space, i.e., it contains only
finitely many orbits of the Borel subgroup $B\times B$ of $G\times G$. 
More generally, one may consider an equivariant embedding $X$ of a
spherical homogeneous $G/H$, and ask whether the closures in $X$ of
the $B$-orbits in $G/H$ are globally $F$-regular.   
The answer is generally negative as some of these
closures have bad singularities, see Ex.6 in \cite{Bri1}. However, the 
question makes sense for the class of multiplicity-free orbit closures
introduced in \cite{Bri1}, since these are normal and Cohen-Macaulay
(see \cite{Bri1} in characteristic $0$, and \cite{Bri2} in arbitrary
characteristic). In fact, the class of multiplicity-free orbit
closures includes the large Schubert varieties in toroidal 
embeddings of $G$, i.e., those which dominate the canonical
compactification of the associated adjoint semisimple group.

\section{Strong $F$-regularity}

In this section $k$ denotes an algebraically closed field of
characteristic $p>0$ and $R$ denotes a commutative $k$-algebra which
is essentially of finite type, i.e., is isomorphic to a localization
of a finitely generated $k$-algebra.  

Composing the $R$-module structure on an $R$-module $M$ with 
the Frobenius map $F : R \rightarrow R$, $r \mapsto r^p$, defines 
a new $R$-module which we denote by $F_* M$. The module defined 
by iterating this procedure $n$ times will be denoted by $F_*^n M$.
In particular, this defines an $R$-module $F_*^n R$ for each 
positive integer $n$ which as an abelian group coincides with 
$R$ but where the $R$-module structure is twisted by the $n$-th
iterated Frobenius morphism $r \mapsto r^{p^n}$. 

When $s \in R$ and $n$ is a positive integer we may define an 
$R$-module map by
$$F_s^n : R \rightarrow F_*^n R,$$
$$ r \mapsto r^{p^n} s.$$
A {\it splitting} of $F_s^n$ is a $R$-module map $\phi :F_*^n R
\rightarrow R$ such that the composed map $\phi \circ F_s^n$
coincides with the identity map on $R$.

\begin{Definition}({\cite{HH}})
The ring $R$ is {\it strongly $F$-regular} if for each $s \in R$, not
contained in a minimal prime of $R$, there exists a positive integer
$n$ such the map $F_s^n$ is split. The affine scheme Spec$(R)$ is said
to be strongly $F$-regular if $R$ is strongly $F$-regular.
\end{Definition}

It is known (see \cite{HH}) that strongly $F$-regular rings are 
reduced, normal, Cohen-Macaulay, and $F$-rational. Moreover, strongly
$F$-regular rings are weakly $F$-regular, i.e. all ideals are tightly
closed. Being strongly $F$-regular is a local condition in the sense 
that $R$ is strongly $F$-regular if and only if all its local rings
are strongly $F$-regular.

\section{Global $F$-regularity}

In this section $X$ will denote a projective variety over an
algebraically closed field $k$ of characteristic $p>0$. (By a variety,
we mean a separated integral scheme of finite type over $k$; in
particular, varieties are irreducible).
When $\Lcal$ is an ample invertible sheaf on $X$ we define the
associated {\it section ring} to be
$$R = R(X,\Lcal) := \bigoplus_{n \in \Zbb} \Gamma(X, \Lcal^n).$$ 
This ring $R$ is a positively graded, finitely generated
$k$-algebra. We may now state

\begin{Definition} (\cite {Smith})
The projective variety $X$ is {\it globally $F$-regular} 
if the ring $R(X,\Lcal)$ is strongly $F$-regular for some
ample invertible sheaf $\Lcal$ on $X$.
\end{Definition}

It can be shown (see Thm.3.10 in \cite{Smith}) that if $X$ is globally 
$F$-regular then the section ring associated to any ample 
invertible sheaf is strongly $F$-regular. Moreover, if $X$ is globally 
$F$-regular then all its local rings are strongly 
$F$-regular. In particular, 

\begin{Corollary}
If the projective variety $X$ is globally $F$-regular then $X$ is
normal and Cohen-Macaulay. Furthermore, the section ring $R(X,\Lcal)$ 
of any ample invertible sheaf $\Lcal$ on $X$ is normal and Cohen-Macaulay. 
\end{Corollary}

\subsection{Frobenius splitting}
\label{F-splitting}

The absolute Frobenius morphism $F : Y \rightarrow Y$ on 
a scheme $Y$ of finite type over $k$, is the map which is 
the identity on the set of points and where the associated map of 
structure sheafs $F^\sharp : \Ocal_Y \rightarrow F_* 
\Ocal_Y$ is the $p$-th power map. Following \cite{MehtaRamanathan}
we say that $Y$ is {\it Frobenius split} if 
$F^\sharp$ splits as a map of $\Ocal_Y$-modules,
i.e. if there exists an $\Ocal_Y$-linear map 
$\phi : F_* \Ocal_Y \rightarrow  \Ocal_Y$ such
that the composed map $ \phi \circ F^\sharp$ is the 
identity. The map $\phi$ is in this case called a
{\it Frobenius splitting} of $Y$. If $\phi$ is a (Frobenius) 
splitting of $Y$ and $Z$ is a closed subscheme 
of $Y$ with associated  ideal sheaf $\mathcal I$, we 
say that $Z$ is {\it compatibly split} (by $\phi$) 
if $\phi(\mathcal I) \subseteq \mathcal I$. In this 
case $\phi$ induces a splitting of $Z$.
It is easily seen that globally $F$-regular 
varieties are Frobenius split but the converse is in 
general not true. 

Let $D$ denote an effective Cartier divisor on $Y$ and 
let $s$ denote the canonical section of the associated 
invertible sheaf $\Ocal_Y(D)$. When $n$ is a positive 
integer we let $F^n(D)$ denote the $\Ocal_Y$-linear
map
$$F^n(D) : \Ocal_Y \rightarrow F_*^n \Ocal_Y(D), $$
$$ t \mapsto t^{p^n} s. $$
We say that $Y$ is {\it stably Frobenius split along $D$} 
if $F^n(D)$ is split, as a map of $\Ocal_Y$-modules, 
for some $n$. In this case $Y$ is Frobenius split 
as well; the {\it induced Frobenius splitting} is 
given by composing the splitting of $F^n(D)$ with
the map $F_* \Ocal_Y \rightarrow 
F_*^n \Ocal_Y(D)$, $t \mapsto 
t^{p^{n-1}}s$. In case $F^1(D)$ is split we simply say 
that $Y$ is {\it Frobenius split along $D$}. 

In the following lemma we will, for later use, collect a number of
standard facts about the concepts introduced above. 

\begin{Lemma}
\label{froblemma}
Let $Y$ be a scheme of finite type over over $k$. 
\begin{enumerate}
\item If $\phi$ is a Frobenius splitting of $Y$ which 
compatibly splits closed subschemes $Z_1$ and $Z_2$ 
then the scheme theoretic intersection $Z_1 \cap Z_2$
is also compatibly split by $\phi$.
\item If $Z$ is compatibly Frobenius split in $Y$ by 
$\phi$ then every irreducible component of $Z$ is 
compatibly split by $\phi$.
\item Assume that $Y$ is Frobenius split along $D$ and that 
the induced splitting compatibly splits a closed subscheme $Z$.
Assume further that none of the irreducible components 
of $Z$ is contained in the support of $D$. Then $Z$ 
is Frobenius split along $D\cap Z$, where $D \cap Z$ denotes 
the restriction of $D$ to $Z$.
\item Let $D' \leq D$ be effective Cartier divisors on $Y$. Then every
(stable) splitting of $Y$ along $D$ induces a (stable) splitting along
$D'$. Moreover, the induced  Frobenius splittings
of $Y$, defined by these two splittings, coincide. 
\item Let $D$ and $D'$ denote effective Cartier divisors on $Y$.
If $Y$ is stably Frobenius split along both $D$ and $D'$
then $Y$ is also stably split along the sum $D + D'$. 
Furthermore, the stable splitting along $D + D'$ may 
be chosen such that the induced Frobenius splitting
of $Y$ coincides with the one induced by the stable splitting along $D'$.  
\item If $Y$ is Frobenius split along $(p-1)D$ for some effective
Cartier divisor $D$, then $D$ (regarded as the zero subscheme of its
canonical section) is compatibly split by the induced Frobenius splitting.
\end{enumerate}
\end{Lemma}

\begin{proof}
For the proof of (1), (2) and (3) see Prop.1.9 in  \cite{Ram}. For the
proof of (4) let $D''$ denote the effective Cartier divisor $D - D'$
and let $s_D$, $s_{D'}$ and $s_{D''}$ denote the canonical sections of
$\Ocal_Y(D)$, $\Ocal_Y(D')$ and $\Ocal_Y(D'')$ respectively. Let 
$$\phi : F^n_* \Ocal_Y(D) \rightarrow \Ocal_Y$$
denote the splitting of the map 
$F^n(D)$
which exists by assumption. Consider now the diagram
$$\xymatrix{\Ocal_Y \ar[rd]_{F^n(D)} \ar[r]^(.4){(F^n)^\sharp} 
& F^n_*\Ocal_Y  \ar[d]^{s_D} \ar[r]^(.5){s_{D'}} 
& F^n_* \Ocal_Y(D') \ar[d]^{\phi'} \ar[ld]_{s_{D''}} \\
& F^n_* \Ocal_Y(D) \ar[r]_(.6){\phi} &  \Ocal_Y \\
}$$
where $\phi'$ is defined such that the diagram is commutative. It
follows that $\phi'$ defines a stable splitting along $D'$ and that it
induces the same  Frobenius splitting 
$F_* \Ocal_Y \rightarrow \Ocal_Y$ as the  splitting $\phi$
along $D$.

Let now $D$ and $D'$ be as described in (5) and let $\phi :
F^n_* \Ocal_Y(D)  \rightarrow  \Ocal_Y$ and 
$\phi' : F^m_* \Ocal_Y(D')  \rightarrow  \Ocal_Y$ denote 
the associated stable splittings. 
By applying the projection formula and the isomorphism
$F^* \Lcal \simeq \Lcal^p$, when $\Lcal$ is any invertible sheaf on
$Y$, we may define the map
$$\eta :F_*^{n+m} \Ocal_Y( D + p^n D') \simeq 
F^m_* (\Ocal_Y(D') \otimes F_*^n \Ocal_Y(D)) 
\rightarrow F_*^m \Ocal_Y(D'),$$
where the latter map is induced by tensoring $\phi$ with 
$ \Ocal_Y(D')$ and applying the functor $F_*^m$.
The composition of $\eta$ with $\phi'$ then defines 
a stable Frobenius splitting along $D + p^n D'$,
and it is easily checked that the induced Frobenius 
splitting coincides with the one induced by $\phi'$.
The statement now follows from (4).

Next we prove (6). Consider an effective Cartier divisor $D$ and a
Frobenius splitting of $Y$ along $(p-1)D$, defined by the morphism
$$ \phi :  F_*  \Ocal_Y((p-1)D) \rightarrow \Ocal_Y.$$
As the statement of (6) is a local condition we may assume that $Y$ is
affine and that $\Ocal_Y(D) \simeq \Ocal_Y$. Let $s$ be the regular
function on $Y$ associated, under the latter isomorphism, to the
canonical section of $\Ocal_Y(D)$. Then $\phi$ is identified with an
$\Ocal_Y$-linear morphism 
$\tilde \phi : F_*  \Ocal_Y \rightarrow \Ocal_Y,$
and the induced Frobenius splitting is defined by 
$$ \eta : F_*  \Ocal_Y \rightarrow \Ocal_Y, ~ ~
 t \mapsto \tilde \phi(ts^{p-1}).$$
In particular, for $t \in \Ocal_Y(Y)$ it follows that 
$\eta(ts) = s \tilde \phi(t) \in (s)$. Hence the zero subscheme of the
canonical section of $D$, whose ideal is generated by $s$, is
compatibly Frobenius split.
\end{proof}

We also record the following result.

\begin{Lemma}
\label{direct}
Let $f:\tilde{X}\rightarrow X$ be a morphism of projective varieties.
Let $\tilde{Y}$ be a closed subvariety of $\tilde{X}$ and put 
$Y := f(\tilde{Y})$. Assume that $\tilde{X}$ is stably Frobenius
split along an ample effective divisor $\tilde{D}$ not containing
$\tilde{Y}$, and that the induced Frobenius splitting of
$\tilde{X}$ compatibly splits $\tilde{Y}$. 
If the map $f^\sharp : \Ocal_X \rightarrow f_*\Ocal_{\tilde{X}}$ is an
isomorphism, then the map 
$\Ocal_Y \rightarrow f_*\Ocal_{\tilde{Y}}$ is an isomorphism as well.
\end{Lemma}

\begin{proof}
It suffices to show that the composition
$\Ocal_X \rightarrow f_*\Ocal_{\tilde{X}} \rightarrow 
f_*\Ocal_{\tilde{Y}}$ is surjective. This will 
follow if the map  
$$\Gamma(X,\Lcal) \rightarrow 
\Gamma(X,\Lcal \otimes f_*\Ocal_{\tilde{Y}})$$
is surjective for any very ample invertible sheaf $\Lcal$ on $X$.
By the projection formula, this amounts to the surjectivity of the
restriction map 
$$\Gamma(\tilde{X},f^*\Lcal) \rightarrow 
\Gamma(\tilde{Y},f^*\Lcal).$$
The latter map is part of a commutative diagram 
$$\xymatrix{\Gamma(\tilde{X}, f^* \Lcal) \ar[d] \ar[r] 
& \Gamma(\tilde{Y}, f^* \Lcal) \ar[d] \\ 
\Gamma(\tilde{X}, F_*^n (\Ocal_{\tilde{X}}(\tilde{D}) \otimes 
f^* \Lcal^{p^n}))  \ar@/_/[u] \ar[r] 
& \Gamma(\tilde{Y}, F_*^n (\Ocal_{\tilde{X}}(\tilde{D}) \otimes 
f^* \Lcal^{p^n})) \ar@/_/[u]  \\
}$$
where the split vertical maps are induced from the stable Frobenius
splitting of $\tilde{X}$ (resp. $\tilde{Y}$) along $\tilde{D}$ (resp.
$\tilde{D} \cap \tilde{Y}$ using Lemma \ref{froblemma}(3)), and where
the horizontal maps are restriction maps. 
By Prop.3 in \cite{MehtaRamanathan} the lower 
horizontal map is surjective. Hence, as the 
splittings 
of the vertical maps are compatible, we conclude
that the upper horizontal map is also surjective.
\end{proof}

Assume now that $Y$ is a nonsingular variety and let $\omega_Y$ denote
its dualizing sheaf. By duality for the finite morphism  $F$ it
follows that
$$ {\mathcal Hom}_{\Ocal_Y}(F_* \Ocal_Y, \Ocal_Y) 
\simeq F_*(\omega_Y^{1-p}),$$ 
which means that a Frobenius splitting of $Y$ is 
the same as a global section of $\omega_Y^{1-p}$
with certain properties. More precisely, let 
$$ C : F_* (\omega_Y^{1-p}) \rightarrow \Ocal_Y,$$
$$ s \mapsto s(1),$$
be the morphism defined by the isomorphism above. Then 
a global section $s$ of $\omega_Y^{1-p}$ defines a 
Frobenius splitting if and only if $C(s)$ coincides 
with the constant function 1 on $Y$. Assume that 
$s$ is a global section of $\omega_Y^{1-p}$ which 
defines a Frobenius splitting, and let $D$ denote the 
divisor of zeroes of $s$. Then, by the discussion 
above, the composed map $C \circ F^1(D)$ is the 
identity map on $\Ocal_Y$ and hence $C$ defines 
a Frobenius splitting of $Y$ along $D$.

\subsection{A criterion for global $F$-regularity}

The following important result by Smith (see Thm.3.10 in 
\cite{Smith}) connects global $F$-regularity, Frobenius splitting 
and strong $F$-regularity.

\begin{Theorem}
\label{Thm-Smith}
If $X$ is a projective variety over $k$ then 
the following are equivalent

\begin{enumerate}
\item $X$ is globally $F$-regular. 

\item $X$ is stably Frobenius split along an ample 
effective Cartier divisor $D$ and the (affine) 
complement $X \setminus D$ is strongly $F$-regular.

\item $X$ is stably Frobenius split along every 
effective Cartier divisor. 
\end{enumerate}
\end{Theorem}
The connection between (1) and (3) in this theorem 
leads to the following result which can be found 
in \cite{LauPedTho}.

\begin{Corollary}
\label{pushforward}
Let $f : \tilde{X} \rightarrow X$ be a morphism of projective 
varieties. If the map 
$f^\sharp : \Ocal_X \rightarrow f_* \Ocal_{\tilde{X}}$ 
is an isomorphism and $\tilde{X}$ is globally $F$-regular then 
$X$ is also globally $F$-regular.
\end{Corollary}

\section{Equivariant embeddings of reductive groups}

In this section $G$ will denote a connected reductive algebraic group
over an algebraically closed field $k$ of arbitrary characteristic. 
We will fix a Borel subgroup $B$ and a maximal torus $T \subseteq B$
of $G$. The Weyl group $N_G(T)/T$ will be denoted by $W$. For any
$w\in W$ we denote by 
$\dot{w}$ a representative in $N_G(T)$. The set of roots defined by
$T$ will be denoted by $\Phi$. To each root $\alpha$ is associated a
reflection $s_\alpha$ in $W$. We choose the set of positive roots
$\Phi^+$ to consist of the roots in $\Phi$ defined by $B$, i.e. $\Phi^+$
consists of the $T$-weights of the Lie algebra of the unipotent  
radical of $B$. The set of positive simple roots will be denoted by
$\Delta$ and the associated simple reflections will be denoted by
$s_1, \dots, s_\ell$. Each element $w$ in $W$ is a product of simple 
reflections and the least number of factors needed in such a product
will be denoted by $\ell(w)$ and will be called the length of $w$. 
The unique element in $W$ of maximal length is denoted by $w_0$. 

We will denote by $\Lambda$ the character group of $T$ and by
$\Lambda^+$ the subset of dominant weights (i.e., those characters
having nonnegative scalar product with all the simple coroots). We
have a partial ordering $\le$ on the group $\Lambda$, where 
$\mu \le \lambda$ if and only if $\lambda -\mu$ is a linear
combination of the simple roots with nonnegative integer coefficients.

For any $w\in W$, the double coset $B \dot{w} B$ is a locally closed
subvariety of $G$ which only depends on $w$; we will denote this 
subvariety by $BwB$. By the Bruhat decomposition the group $G$ is the
disjoint union of the double cosets $BwB$, $w\in W$. Moreover, 
$\dim(BwB) = \ell(w) + \dim(B) = \ell(w) + \ell(w_0) + \ell$. The
closure in $G$ of any $BwB$ is the union of the $BvB$, where $v\in W$
and $v\le w$ for the Bruhat ordering of $W$. 

An {\it equivariant embedding} of $G$ is a normal 
$G \times G$-variety $X$ containing $G$ as an open 
subset and where the induced $G \times G$-action 
on $G$ is given by left and right translation.
(In other words, $X$ is an equivariant embedding
of the homogeneous space $G\times G/\diag G \simeq G$.)

The {\it boundary} of the equivariant embedding $X$ is the closed
$G\times G$-stable subset $X \setminus G$, denoted by 
$\partial X$. Its irreducible components $D_1,\dots, D_n$ are the 
{\it boundary divisors}; they are indeed of codimension $1$, as the
open subset $G$ is affine (see Prop.II.3.1 in \cite{Hartshorne}).

When $X$ is an equivariant embedding of $G$ we denote by $X(w)$, 
$w\in W$, the closure in $X$ of the double coset $BwB$ (in particular,
$X(w_0)=X$). In this section we want to study the geometry of these
{\it large Schubert varieties} $X(w)$. Those of codimension $1$ are
the $X(w_0s_i)$, $i=1, \dots, \ell$; they will be denoted by 
$X_1, \dots, X_l$.

\subsection{Two preliminary geometric results}

A key ingredient in our study is the following

\begin{Proposition}
\label{ample}
Let $X$ be a projective embedding of $G$, with boundary divisors
$D_1,\ldots,D_n$. 

\begin{enumerate}
\item There exists a very ample $G\times G$-linearized invertible
sheaf $\Lcal$ over $X$ such that $\Spec R(X,\Lcal)$ is an affine
embedding of the group $G \times \Gbb_m$, where action of the
multiplicative group $\Gbb_m$ on $\Spec R(X,\Lcal)$ corresponds to the
grading of $R(X,\Lcal)$.

\item There exist positive integers $a_1,\dots,a_n$ such that the
invertible sheaf $\Ocal_X(\sum_{i=1}^n a_iD_i)$ is ample. 
\end{enumerate}
\end{Proposition}

\begin{proof}
(1) We may find a very ample $G\times G$-linearized invertible sheaf
$\Lcal$ on $X$, see e.g. Cor.1.6 in \cite{Mum}.
Then the pull-back of $\Lcal$ to the open orbit
$G \simeq G\times G/\diag G$ is the linearized invertible sheaf
associated with a character of the isotropy group $\diag G$. Such a
character extends to a character of $G\times G$, so that (changing the
linearization) we may assume that the pull-back of $\Lcal$ to $G$ is
trivial as a linearized invertible sheaf.

Replacing $\Lcal$ with some positive power, we may also assume that
the ring $R(X,\Lcal)$ is normal. Then $\hat X := \Spec R(X,\Lcal)$ is
a normal affine variety endowed with an action of 
$G\times G \times \Gbb_m$, where $\Gbb_m$ acts via the grading of
$R(X,\Lcal)$. By our assumptions on $\Lcal$, the affine cone $\hat X$
is an affine embedding of the group $G\times \Gbb_m=:\hat G$.

(2) Let $\hat X$, $\hat G$ as above. Then $\hat X$ is a linear
algebraic monoid with unit group $\hat G$, by Prop.1 in \cite{Rit1}. 
So, by Thm.3.15 in \cite{Put}, $\hat X$ admits an embedding into some
matrix ring $M_n(k)$ as a closed submonoid (with respect to the
multiplication of matrices). We claim that $\hat G$ identifies with
$\hat X \cap \GL_n(k)$ under this embedding. Indeed, the inclusion
$\hat G \subseteq \hat X \cap \GL_n(k)$ is clear. Conversely,
if $\gamma \in \hat X \cap \GL_n(k)$ then the images 
$\gamma^i \hat X$ form a decreasing sequence of closed subsets of 
$\hat X$. Thus, $\gamma^i \hat X = \gamma^{i+1}\hat X$ for $i\gg 0$. 
It follows that $\hat X= \gamma \hat X$, whence $\gamma$ has a right
inverse. Likewise, $\gamma$ has a left inverse, which completes the
proof of the claim.

Let $s$ be the regular function on $\hat X$ given by the restriction
of the determinant function on $M_n(k)$. By the claim, the zero set of
$s$ is precisely the boundary $\partial \hat X$. Further, $s$ is an
eigenvector of $\Gbb_m$, by the multiplicative property of the
determinant. So $s$ is a section of a positive power of 
$\Lcal$, with zero set being $\partial X$.
\end{proof}

We also recall the following result which is known under a stronger
form (Prop.3 in \cite{Rit2}, see also Prop.6.2.5 in \cite{BriKum}).

\begin{Lemma}
\label{resolution}
For any equivariant embedding $X$ of $G$, there exists a nonsingular
equivariant embedding $\tilde{X}$ and a projective morphism 
$$ f : \tilde{X} \rightarrow X$$
which induces the identity on $G$.
\end{Lemma}

We will refer to $f$ as an {\it equivariant resolution} of $X$. Note
that $f$ is $G\times G$-equivariant and birational. Since $X$ is
assumed to be normal, it follows that
$f_*\Ocal_{\tilde{X}} =\Ocal_X$. Further, $\tilde{X}$ is
projective if $X$ is. Also, note that $f$ restricts to a birational
morphism $\tilde{X}(w) \rightarrow X(w)$, for any $w\in W$. 
Together with Lemma \ref{direct}, this will allow us to reduce
questions on $X(w)$ to the case where $X$ is nonsingular.

\subsection{Frobenius splitting of nonsingular embeddings}

Now fix a nonsingular equivariant embedding $X$ of $G$. We assume from 
now on that the ground field $k$ has characteristic $p>0$. It is known
(see \cite{Rit2}, or Prop.6.2.6 in \cite{BriKum}) that the inverse of
the dualizing sheaf on $X$ equals 
$$\omega_X^{-1} \simeq 
\Ocal_X(\partial X + \sum_{i=1}^l (X_i + \tilde{X_i})) ,$$
where $\tilde X_i = (\dot{w_0},\dot{w_0}) X_i$, and that 
the $(p-1)$-th power $s_X^{p-1}$ of the canonical section $s_X$
of the right hand side defines a Frobenius splitting of $X$. 
As noticed at the end of Section \ref{F-splitting}, this yields in
fact a splitting of $X$ along $D$, where
$$D := (p-1)(\partial X + \sum_{i=1}^l (X_i + \tilde{X_i})).$$
This leads to the following result.

\begin{Proposition}
\label{D-splitting}
Let $X$ denote a nonsingular equivariant embedding of $G$ over a field 
of characteristic $p>0$. Then $X$ is Frobenius split along
$(p-1)\partial X$, compatibly with the large Schubert subvarieties
$X(w)$, $w\in W$.
\end{Proposition}

\begin{proof}
Denote by $\eta : F_* \Ocal_X \rightarrow \Ocal_X$ the underlying
Frobenius splitting of $X$ induced by $s_X$. By Lemma 
\ref{froblemma}(2)(6), each $X_i$ is compatibly Frobenius 
split by $\eta$.
In other words, $\eta$ is compatible with the $X(w)$, where 
$\ell(w) = \ell(w_0)-1$. Now consider $w\in W$ such that 
$\ell(w) \le \ell(w_0)-2$. By Lem.10.3 in \cite{BGG}, there exist
distinct $w_1$, $w_2$ in $W$ such that $w<w_1$, $w<w_2$, and 
$\ell(w_1)=\ell(w_2)=\ell(w)+1$. Then $X(w)$ is contained in
$X(w_1)\cap X(w_2)$ as an irreducible component. Now Lemma
\ref{froblemma}(1)(2) implies by decreasing induction on $\ell(w)$ 
that $X(w)$ is compatibly Frobenius split by $\eta$. That 
$\eta$ is induced by a Frobenius splitting of $X$ along 
$(p-1)\partial X$ follows from Lemma \ref{froblemma}(4).
\end{proof}

\subsection{The main results}

We still assume that $k$ has characteristic $p>0$.

\begin{Theorem}
\label{main}
Let $X$ denote a projective equivariant embedding  of $G$.
Then each $X(w)$, $w \in W$, is globally $F$-regular. 
\end{Theorem}

\begin{proof}
First we consider the case where $X$ is nonsingular. Then, 
by Lemma \ref{froblemma}(3) and Proposition \ref{D-splitting} each
$X(w)$ is Frobenius split along $(p-1)(\partial X \cap X(w))$. 
Together with Lemma \ref{froblemma}(4)(5), it follows that $X(w)$ is
stably Frobenius split along any divisor 
$\sum_{i=1}^n a_i (D_i\cap X(w))$, with $a_i >0$. By Proposition
\ref{ample}, we may find such a divisor which is ample on $X$. Then
the restriction 
$$ D =   \sum_{i=1}^n a_i (D_i \cap X(w)), $$
is an effective ample Cartier divisor on $X(w)$ with support 
$\partial X \cap X(w)$. Hence, by Theorem \ref{Thm-Smith} it is enough
to prove that the open affine subset 
$$G(w) = X(w) \setminus \partial X$$
is strongly $F$-regular. 

Notice that the set $G(w)$ coincides with the closure of $Bw B$
in $G$. Hence, there is a surjective map  
$$\pi(w) : G(w) \rightarrow S(w) \subseteq G/B, $$ 
onto the corresponding Schubert variety $S(w)$. By
\cite{LauPedTho} $S(w)$ is globally $F$-regular and hence locally
strongly $F$-regular. Moreover, by the Bruhat decomposition there
exists a covering of $S(w)$ by open affine 
subsets $U_i$, $i \in I$, such that $\pi(w)^{-1}(U_i) 
\simeq U_i \times B$. As $B$ is  smooth
and  $U_i$ is strongly $F$-regular it follows 
that $U_i \times B$ is strongly $F$-regular 
(Lem.4.1 in \cite{LyuSmith}). Hence, the affine variety 
$G(w)$ is also strongly $F$-regular. This completes the proof
in the case of nonsingular $X$.

In the general case, we may choose an equivariant resolution 
$f : \tilde{X} \rightarrow X$ (Lemma \ref{resolution}). By 
the considerations above the equivariant embedding $\tilde{X}$ 
is stably Frobenius split along an ample Cartier divisor. 
Furthermore according to the last part of Lemma \ref{froblemma}(4), 
this stable splitting may be chosen such that each $\tilde{X}(w)$ 
is compatibly Frobenius split. Then, by Lemma \ref{direct}, the map 
$\Ocal_{X(w)} \rightarrow f_*\Ocal_{\tilde{X}(w)}$
is an isomorphism and the global $F$-regularity of $X(w)$ 
hence follows from Corollary \ref{pushforward}.
\end{proof}

\begin{Corollary}
\label{affp}
Let $X$ denote an affine equivariant embedding of $G$. Then each
$X(w)$, $w\in W$, is strongly $F$-regular.
\end{Corollary}

\begin{proof}
We may embed $X$ as a closed $G\times G$-stable subvariety of a
$G\times G$-module $M$. Let $\overline{X}$ be the normalization of the
closure of $X$ in the projectivization of $M\oplus k$. Then
$\overline{X}$ is a projective equivariant embedding of $G$
containing $X$ as an open affine subset. By Theorem \ref{main}
each $\overline{X}(w)$ is globally $F$-regular. In particular,
every local ring of $\overline{X}(w)$ is strongly $F$-regular.
As $X(w)$ is an open subset of $\overline{X}(w)$ this implies 
that every local ring of the affine variety $X(w)$ is 
strongly $F$-regular. This proves the claim as the condition 
of being strongly $F$-regular is local.
\end{proof}

\begin{Corollary}
\label{singp}
Let  $X$ denote any equivariant embedding of $G$. Then each $X(w)$, 
$w\in W$, is normal and Cohen-Macaulay.
\end{Corollary}  

\begin{proof}
This follows from Theorem \ref{main} by using that $X$ has an open
cover by equivariant embeddings which are also open subsets of
projective equivariant embeddings, see \cite{Sum1,Sum2}.
\end{proof}

\subsection{From characteristic $p$ to characteristic $0$}

In this section, $k$ is of characteristic $0$. We will obtain versions
of Theorem \ref{main} and of Corollaries \ref{affp}, \ref{singp}, by
using the notions of strongly (resp.~globally) $F$-regular type
(\cite{Smith}) that we briefly review. 

Let $Y$ be a scheme of finite type over $k$. Then $Y$ is defined over
some finitely generated subring $A$ of $k$. This yields a scheme $Y_A$
which is flat and of finite type over $\Spec(A)$, such that $Y$ is 
naturally identified with the scheme $Y_A\times_{\Spec(A)} \Spec(k)$. 
On the other hand,
the geometric fibers of $Y_A$ at closed points of $\Spec(A)$ are
schemes over algebraic closures of finite fields (of various
characteristics).

\begin{Definition} (\cite{Smith}) 
An affine (resp.~projective) variety $X$ is 
{\it of strongly (resp.~globally) $F$-regular type} if  
$X$ is defined over some finitely generated subring $A$ of $k$
such that the geometric fibers of $X_A$ over a dense subset of
closed points of $\Spec(A)$ are strongly (resp.~globally) $F$-regular.
\end{Definition}

Remember that any strongly (resp. globally) $F$-regular 
variety is locally $F$-rational. It follows that any 
variety $X$ of strongly (resp.~globally) $F$-regular type is 
of $F$-rational type (this latter notion is defined 
similarly to the definition of strongly/globally $F$-regular
type). Hence by Thm4.3. in \cite{Smith2} it follows 
that $X$ has rational singularities, in particular, 
$X$ is normal and Cohen-Macaulay. 

\begin{Theorem}
Let $X$ denote an affine (resp.~projective) equivariant embedding of
$G$ over a field of characteristic $0$. Then any $X(w)$, $w\in W$, is
of strongly (resp.~globally) $F$-regular type.
\end{Theorem}

\begin{proof}
By Proposition \ref{ample} (1), it suffices to treat the affine
case. For this, we will recall the classification of affine
equivariant embeddings, after \cite{Vin} (generalized in \cite{Rit1}
to arbitrary characteristic), and show that any such embedding $X$ is
defined and flat over $\Spec(\Zbb)$.

Put $R:=\Gamma(G,\Ocal_G)$ and $S:=\Gamma(X,\Ocal_X)$, then $S$ is a
$G\times G$-stable subalgebra of $R$. Further, $S$ is finitely
generated and normal, with the same quotient field as $R$.
Recall the isomorphism of $G\times G$-modules
$$
R \cong \bigoplus_{\lambda\in\Lambda^+} 
\nabla(\lambda) \otimes \nabla(-w_0\lambda),
$$
where $\nabla(\lambda)$ denotes the simple $G$-module with highest
weight $\lambda$. It follows that 
$$
S \cong \bigoplus_{\lambda\in\Mcal} 
\nabla(\lambda) \otimes \nabla(-w_0\lambda),
$$
for some subset $\Mcal$ of $\Lambda^+$. Thus, the weights of 
$T\times T$ in the invariant subring $S^{U\times U}$ are exactly the
$(\lambda,-w_0\lambda)$, where $\lambda\in\Mcal$; each such weight
has multiplicity $1$. Since $S^{U\times U}$ is a finitely generated,
normal domain (see e.g. \cite{Gro}), the corresponding affine variety
is a toric variety for the left $T$-action. Thus, $\Mcal$ 
is the intersection of $\Lambda$ with a rational polyhedral convex
cone of nonempty interior in $\Lambda\otimes_{\Zbb}\Rbb$, contained in
the positive chamber. 

One may show that $\Mcal$ satisfies the following saturation property:
For any $\lambda\in\Mcal$ and $\mu\in\Lambda^+$ such that
$\mu\le\lambda$, then $\mu\in\Mcal$. Conversely, any $\Mcal$
satisfying the preceding properties yields an affine embedding of $G$, 
see \cite{Vin}.

Next let $G_{\Zbb}$ be the split $\Zbb$-form of $G$, with affine
coordinate ring $R_{\Zbb}$. For any ring $A$, this defines the ring
$R_A:=R_{\Zbb}\otimes_{\Zbb} A$ and the corresponding group $G_A$.
In particular, we obtain the $\Qbb$-form $R_{\Qbb}$ of $R$. Now the
preceding decomposition of $R$ is defined over $\Qbb$; further, the
subspace
$$
S_{\Qbb} := S\cap R_{\Qbb} = \bigoplus_{\lambda\in\Mcal} 
\nabla_{\Qbb}(\lambda) \otimes \nabla_{\Qbb}(-w_0\lambda)
$$
(with obvious notation) is a subalgebra of $R_{\Qbb}$, and a
$\Qbb$-form of $S$. Put $S_{\Zbb}:=S_{\Qbb} \cap R_{\Zbb}$ (then the
quotient $R_{\Zbb}/S_{\Zbb}$ is torsion-free), and 
$$
R_p:=R_{\Zbb}\otimes_{\Zbb}\overline{\Fbb_p}, \quad
S_p:=S_{\Zbb}\otimes_{\Zbb}\overline{\Fbb_p},
$$
where $p$ is any prime number, $\Fbb_p$ denotes the field with $p$
elements, and $\overline{\Fbb_p}$ denotes its algebraic closure.
Define likewise the connected reductive group $G_p$ over
$\overline{\Fbb_p}$ and its subgroups $B_p$, $T_p$, $U_p$. Then
$R_p=\Gamma(G_p,\Ocal_{G_p})$, and $S_p$ is a
$\overline{\Fbb_p}$-subalgebra of $R_p$, stable under the action of
$G_p\times G_p$. We will show that $S_p$ is the coordinate ring of an
affine equivariant embedding $X_p$ of $G_p$.

By Prop.II.4.20 in \cite{Jantzen} (see also Thm.4.2.5 in
\cite{BriKum}), the $G_p\times G_p$-module 
$R_p$ has an increasing filtration with subquotients being the 
$\nabla_p(\lambda)\otimes \nabla_p(-w_0\lambda)$
($\lambda\in\Lambda^+)$, where now $\nabla_p(\lambda)$ denotes 
the dual Weyl module of highest weight $\lambda$. Further, the 
proof of this result given in \cite{BriKum} also shows that the 
$G_p\times G_p$-module $S_p$ has an increasing filtration with
subquotients being the $\nabla_p(\lambda)\otimes \nabla_p(-w_0\lambda)$
($\lambda\in\Mcal)$. In particular, this module has a good
filtration. Using Lem.II.2.13 and Prop.II.4.16 in \cite{Jantzen}, it
follows that the weights of $T_p\times T_p$ in the invariant subring
$S_p^{U_p\times U_p}$ are again the $(\lambda,-w_0\lambda)$, where
$\lambda\in\Mcal$; each such weight has multiplicity $1$. Therefore,
the algebra $S_p^{U_p\times U_p}$ is finitely generated and normal. By
\cite{Gro}, the algebra $S_p$ is finitely generated and normal as well. 

Put $X_p:=\Spec(S_p)$, then $X_p$ is a normal affine variety where
$G_p\times G_p$ acts with a dense orbit. We now show that this orbit
is isomorphic to $G_p\times G_p/\diag G_p$; equivalently, the
morphism $G_p \to X_p$ associated with the inclusion 
$S_p\subseteq R_p$ is an open immersion. Since the corresponding
morphism $G\to X$ is an open immersion, we may find 
$f\in S^{U\times U}$ with zero set the complement of the open 
$B\times B$-orbit $Bw_0B$. Replacing $f$ with a scalar multiple, we
may assume that $f\in S_{\Zbb}$ is a lift of a nonzero 
$f_p\in S_p^{U_p\times U_p}$. Then 
$R[f^{-1}] = S[f^{-1}] = \Gamma(Bw_0B,\Ocal_{Bw_0B})$, so that 
$R_{\Zbb}[f^{-1}] = S_{\Zbb}[f^{-1}]$. Thus, 
$R_p[f_p^{-1}] = S_p[f_p^{-1}]$; equivalently, the morphism 
$(f_p\ne 0) = B_p w_0 B_p \to X_p$ is an open immersion. Since the 
$G_p\times G_p$-translates of $B_p w_0 B_p$ cover $G_p$, we have shown
that the reduction $X_p$ is an equivariant embedding of
$G_p$. Further, since all the double classes $BwB$ in $G$ are defined
over $\Zbb$, their closures $X(w)$ in $X$ are also defined over
$\Zbb$, with reductions $X(w)_p$.
\end{proof}

By the argument of Corollary \ref{singp}, this implies readily

\begin{Corollary}
\label{sing0}
Let $X$ denote an equivariant embedding of $G$ over a field of
characteristic $0$. Then each $X(w)$, $w \in W$, has rational
singularities.
\end{Corollary}

\bibliographystyle{amsplain}

\end{document}